\documentclass[12pt]{article}
\usepackage{amsmath}
\usepackage{amssymb}
\usepackage{graphicx}

\begin{document}

\newcommand \be  {\begin{equation}}
\newcommand \bea {\begin{eqnarray} \nonumber }
\newcommand \ee  {\end{equation}}
\newcommand \eea {\end{eqnarray}}

\title{  Counting stationary points of the loss function in the simplest
constrained least-square optimization \footnote{The text is based on
the presentation at the workshop ``Random Matrix Theory: Applications in the Information Era'', 29 Apr 2019 -- 3 May 2019  2019, Krakow, Poland }}

\vskip 0.2cm
\author{Yan V. Fyodorov $^1$ and Rashel Tublin $^1$\\
\noindent\small{$^1$ King's College London, Department of Mathematics, London  WC2R 2LS, United Kingdom}}

\maketitle

\begin{abstract}
We use Kac-Rice method to analyze statistical features of an ``optimization landscape'' of the loss function in a random  version of the Oblique Procrustes Problem, one of the simplest optimization problems of the least-square type on  a sphere.

\end{abstract}

\section{Introduction}
One of the simplest optimization problems of the  least-square type arising in the Multiple Factor Data Analysis  is the following \\ [0.5ex]
{\sf Oblique Procrustes Problem} \cite{Browne1967}: {\it  For a given pair of
$M\times N$ matrices $\mathbf{A}$ and $\mathbf{B}$ find such $N\times N$ matrix $\mathbf{ X}$ that the equality $\mathbf{B=AX}$  holds as close as possible and columns ${\bf x}_i\in \mathbb{R}^N, \, i=1,\ldots N$ are all of the same fixed length: $||{\bf x}||_2:=\sqrt{\sum_i x_i^2}=const$. }\\[0.5ex]
For $M>N$ the associated system of linear equations is over-complete
and a solution can be found separately for each column  ${\bf x}$ by minimizing
the loss/cost function
\be\label{lossfun}
H({\bf x})=\frac{1}{2}||A{\bf x}-{\bf b}||^2:=\frac{1}{2}\sum_{k=1}^M\left[\sum_{j=1}^NA_{kj}{\bf x}_j-b_k\right]^2, \quad ||{\bf x}||_2:=const\,.
\ee

The problem was first analysed in that setting by Browne in 1967 \cite{Browne1967}, and then independently by numerical mathematicians (see e.g.  \cite{Gander1981,GolubMatt1991}) who used the Lagrange multiplier to take care of the spherical constraint. Introducing the Lagrangian ${\cal L}_{\lambda,{\bf s}}({\bf x})={\cal H}({\bf x})-\frac{\lambda}{2}({\bf x},{\bf x})$, with real $\lambda$ being the Lagrange multiplier, the stationary conditions $\nabla {\cal L}_{\lambda,{\bf s}}({\bf x})=0$ yield a linear system:
\be\label{minLagr}
 A^T\left[A{\bf x}-{\bf b}\right]=\lambda {\bf x}, \quad \Rightarrow {\bf x}=(A^TA-\lambda I_N)^{-1} A^T{\bf b}
\ee
We find it convenient to use the normalization such that the radius of the
sphere is $ ||{\bf x}||_2:=\sqrt{N}$, with the spherical constraint  yielding the equation for $\lambda$ in the form:
\be\label{Lagrange}
{\bf b}^TA\,\frac{1}{\left(A^TA-\lambda I_N\right)^{2}}\,A^T{\bf b}=N
\ee
which is equivalent to a  polynomial equation of degree $2N$ in $\lambda$. Each  real solution for the  Lagrange multiplier $\lambda_i$ corresponds to a stationary point ${\bf x}_i$ of the loss function $H({\bf x})=\frac{1}{2}||A{\bf x}-{\bf b}||^2$ on the sphere ${\bf x}^2=N$ and one can show that the order $\lambda_1<\lambda_2<\ldots <\lambda_{\cal N}$ implies $H({\bf x}_1)<H({\bf x}_2)<\ldots <H\left({\bf x}_{\cal N}\right)$ \cite{Browne1967}. Thus the  minimal loss is given by ${\cal E}_{min}=H\left({\bf x}_1\right)$.

Actually, the loss function (\ref{lossfun}) is one of the simplest examples of the ``optimization landscape'', interest in which governs developing various search algorithms efficiently converging to the global minimum.
To consider a ``typical'' landscape it makes sense to assume that the parameters of the model, i.e. the  matrix $A$ and the vector $b$, are random.
Geometrical and topological properties of random landscapes have general and intrinsic mathematical interest, see e.g. \cite{FLL2015}, and attracted considerable attention in recent years due to their relevance in the area of ``deep learning`` and optimization, see e.g. \cite{lossmultilayer,RBBC2018}.
Fruitful analogies with spin glasses where ``energy landscapes'' have been under intensive investigation for some time, see  \cite{Auf1,AufBenCer13,BraDea07,Fyo04,Fyo15,Fyo16,FyoNad12,FyoWil07}, plays an important role in guiding the intuition in this area.
In this context the goal of the present research is to investigate the simplest landscape Eq.(\ref{lossfun}) by counting the stationary points  via  the Lagrange multipliers $\lambda_i, \, i=1,\ldots,{\cal N}\le 2N$ and eventually find the  minimal loss ${\cal E}_{min}$.  For concreteness and analytical tractability we assume the entries $A_{kj}$ of $M\times N, \, M>N$ matrix $A$ to be i.i.d. normal real variables such that
$A^TA=W$ is $N\times N$   Wishart with the probability density
\be\label{Wishart}
 P_{N,M}(W)=C_{N,M}e^{-\frac{N}{2}\mbox{\small Tr} W}\left( \det W\right)^{\frac{M-N-1}{2}}
\ee
 We will also assume for convenience that  the vector ${\bf b}$ is
normally distributed: ${\bf b}=\sigma \, \mathbf{\xi}$  with $\sigma>0$ and the components of  $\mathbf{\xi}=(\xi_1, \dots , \xi_M)^T$ are mean zero standard normals.

\section{Qualitative considerations and the  Kac-Rice method.}
The equation Eq.(\ref{Lagrange}) for the Lagrange multiplier can be conveniently written
in terms of  $N$ nonzero eigenvalues $s_{1}, \ldots,s_{N}$ of $M\times M$  matrix  $W^{(a)}=AA^T$
and the associated eigenvectors ${\bf v}_{i}$:\\[1ex]
\be\label{equlambda_new_xi_1}
    \sum\limits_{i=1}^N \frac{s_i}{(\lambda-s_i)^2} (\mathbf{\xi}^T \mathbf v_i)^2= \frac{N}{\sigma^2}
\ee
The left-hand side is a positive function of $\lambda$ having a single minimum between every consecutive pair of eigenvalues of $W^{(a)}$. This implies there are 0 or 2 solutions of \eqref{equlambda_new_xi_1} (and 1 solution with probability zero) for $\lambda$ between every consecutive pair of eigenvalues, plus two more solutions: one in $\lambda\in (-\infty, s_1)$ and another one in $\lambda\in ( s_N, \infty)$. Note that the latter two solutions exist for any value of $\sigma\in [0,\infty]$, whereas by changing $\sigma$ one changes the number of solutions available between consecutive eigenvalues. In particular, in the limit of vanishing noise (i.e. $\sigma\to 0$ hence $||{\bf b}||_2=0$) every stationary point solution for the Lagrange multiplier corresponds to an eigenvalue $s_n$ of the Wishart matrix, with,  ${\bf x}=\pm {\bf e}_n$ being the associated eigenvectors (hence there are $2N$ stationary points).
On the other hand when $\sigma \to \infty$ the ratio $N/\sigma^2$ in the right-hand side becomes smaller than the global minimum of the left-hand side in $[s_1,s_N]$. Then  only two stationary points remain outside that interval. Obviously, in every particular realization the number of stationary points will gradually change between the two limits as a function of growing $\sigma$, forming a staircase ${\cal N}_{st}(\sigma)$.
Let us illustrate this on a simple example in the case of small $N=5$, see Fig. 1.

    \begin{figure}[h!]
\includegraphics[width=.6\textwidth]{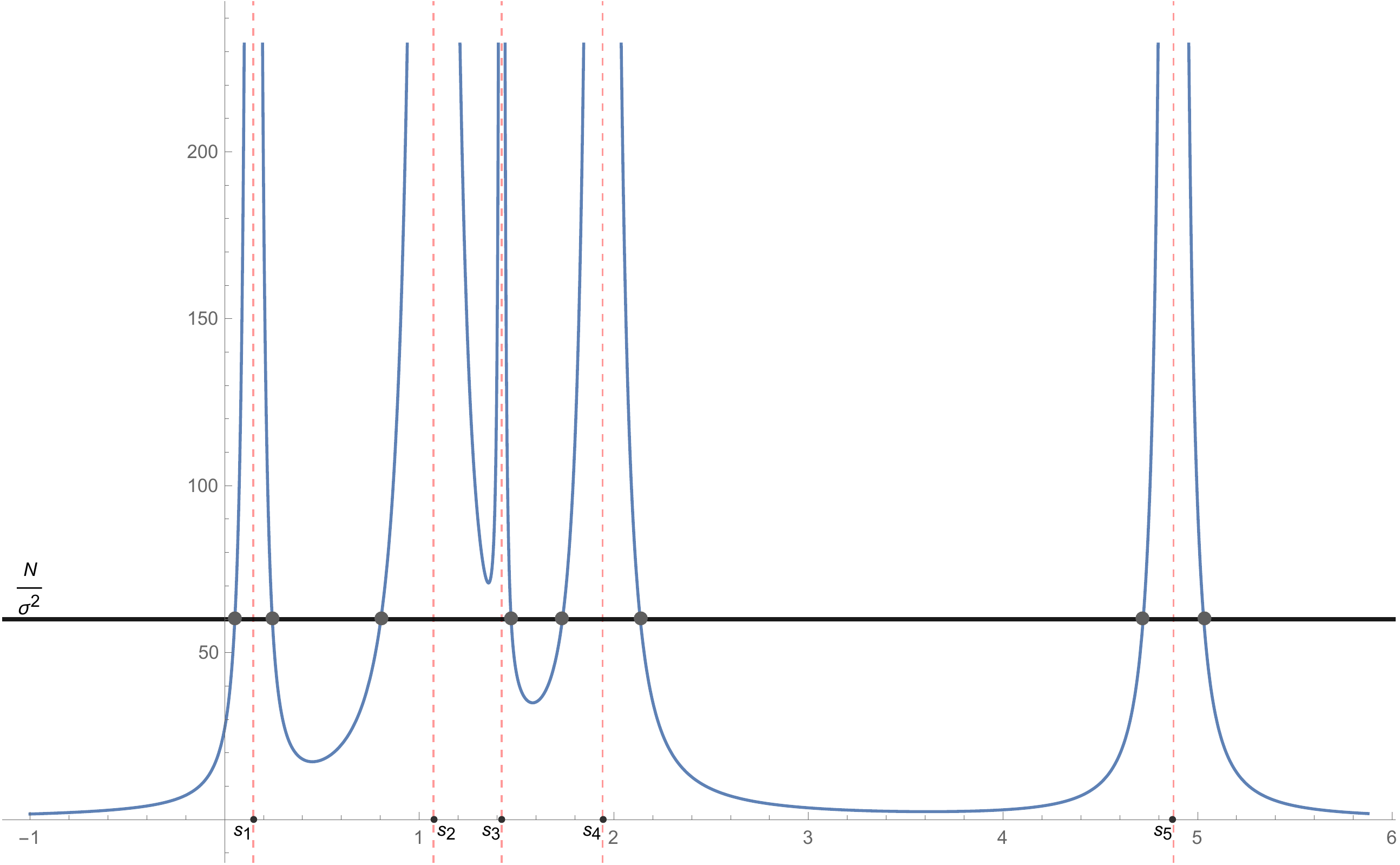}
\includegraphics[width=.4\textwidth]{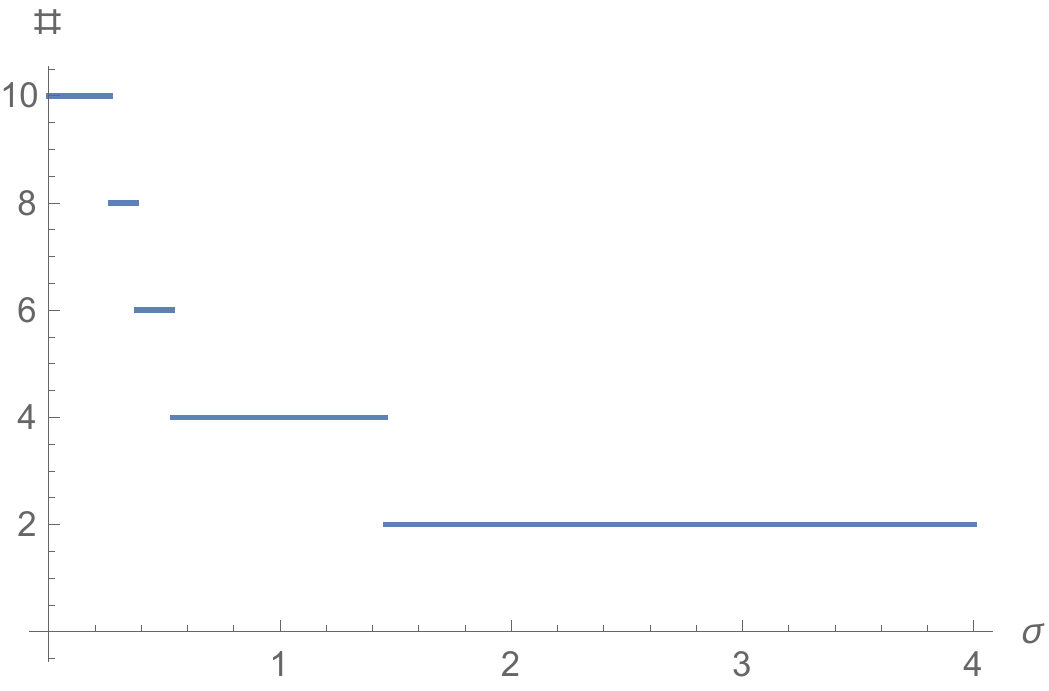}
\label{fig g for N=5}
\centerline{Fig 1. Graphical representation of Eq.\eqref{equlambda_new_xi_1} for $N=5$. }
\end{figure}

This is exactly the ``gradual topology trivialization'' phenomenon
discussed (as the function of magnetic field) for the standard GOE-based spherical model in \cite{FLD2013} (see also \cite{Fyo15,Fyo16})
by adopting formulas derived in the general case by Auffinger et. al. \cite{Auf1,AufBenCer13}.   It is quite easy to see from (\ref{equlambda_new_xi_1}) that the trivialization happens on the scale $\sigma^2\sim 1/N$ as only for such values the left hand-side
is of the same order as the right-hand side for a generic $\lambda \in [s_l,s_{l+1}]$ ( in our normalization the typical
distance $|s_{l}-s_{l+1}|=O(1/N)$). When averaged over the realizations the staircase is replaced by smoothly decreasing function $\left\langle {\cal N}_{st}(\sigma)\right\rangle$ which we will find explicitly using the Kac-Rice approach, and investigate its asymptotics as $N\to \infty$.

The number ${\cal N}_{st}[a,b]$ of real solutions of the Lagrange equation Eq.(\ref{minLagr}), i.e. ${A^T\left[A{\bf x}-{\bf b}\right]-\lambda {\bf x}=0}$  on the sphere ${\bf x}^2=N$ such that $\lambda\in [a,b]$ can be counted by employing the Kac-Rice type formula
\be
     {\cal N}_{st}[a,b]  =  \int_a^b d\lambda \int \delta \left[ A^T \left( A {\bf x} - {\bf b}\right)-\lambda {\bf x} \right] \delta \left({\bf x}^2 - N  \right)
\ee
\[   \times    \left| \det \left(\begin{array}{cc}
            A^TA-\lambda I_N & {\bf x}\\
            -2{\bf x}^T & 0
        \end{array}\right)
    \right| d{\bf x}
\]

 Using Gaussianity of both the matrix entries  $A_{ij}\sim {\cal N}(0,1)$ and the vector components ${\bf b}\sim {\cal N}_M(0,I_M\sigma^2)$ and introducing the parameter $\delta=\frac{1}{2}\ln{(1+\sigma^2)}$ one can eventually find the mean number of solutions as $$\mathbb{E}\left\{{\cal N}_{st}[a,b]\right\}=\int_a^b p(\lambda)\,d\lambda$$
 with the density $p(\lambda)$ for $\lambda>0$  given by
\be p(\lambda\ge 0)= 2 \sqrt{\frac{N}{\pi}} \frac{e^{-\frac{M+N-1}{2}\delta}}{\sqrt{\sinh{\delta}}}  K_{\frac{M-N}{2}} \!\left( \frac{N \lambda}{2 \sinh{\delta}}\right) e^{\frac{N\lambda}{2} \coth \delta }
 \left\langle\rho_N(\lambda)\right\rangle \sqrt{\lambda}
 \ee
 where $K_{\nu}(z)$ is the Bessel-Macdonald function, and $\left\langle\rho_N(\lambda)\right\rangle $ stands for the mean eigenvalue density of $N\times N$ real Wishart matrices $W$ distributed according Eq.(\ref{Wishart}). Such density  for any values $M\ge N$  can be found in  \cite{LNVbook}.
For negative values of the Lagrange multiplier $\lambda$ we have instead:
\be p(\lambda<0)=\frac{N!N^{(M-N)/2}}{2^{(M+N-3)/2}} \frac 1{\Gamma \left( \frac N 2  \right) \Gamma \left( \frac M 2  \right)}\frac {e^{-(M+N-1)\delta/2}}{\sqrt{\sinh \delta} }  e^{-\frac 12 N|\lambda|(\coth\delta -1)} |\lambda|^{(M-N)/2}
\ee
\[
\times  \left[ \sum\limits_{j=0}^{N-1} \binom{M-1}{N-1-j} \frac 1{j!} (N|\lambda| )^j  \right] K_{\frac{M-N}2} \left( \frac {N|\lambda|}{2\sinh \delta}   \right)
\]

These formulas are exact, and we can compare them with the direct numerical simulations in Fig.2 for
moderate matrix sizes.

\begin{figure}[h!]
\includegraphics[width=.30\textwidth]{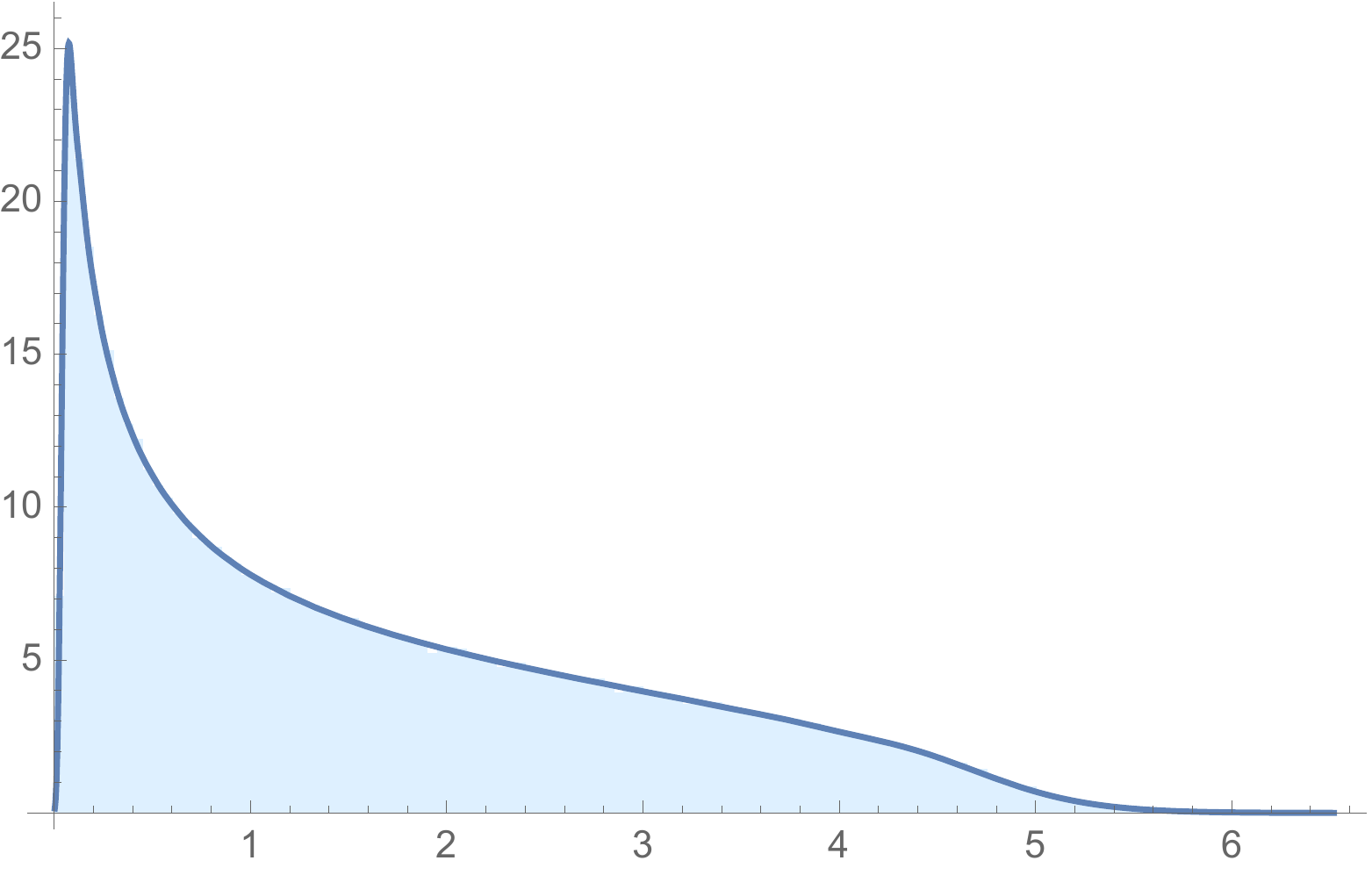}
\includegraphics[width=.30\textwidth]{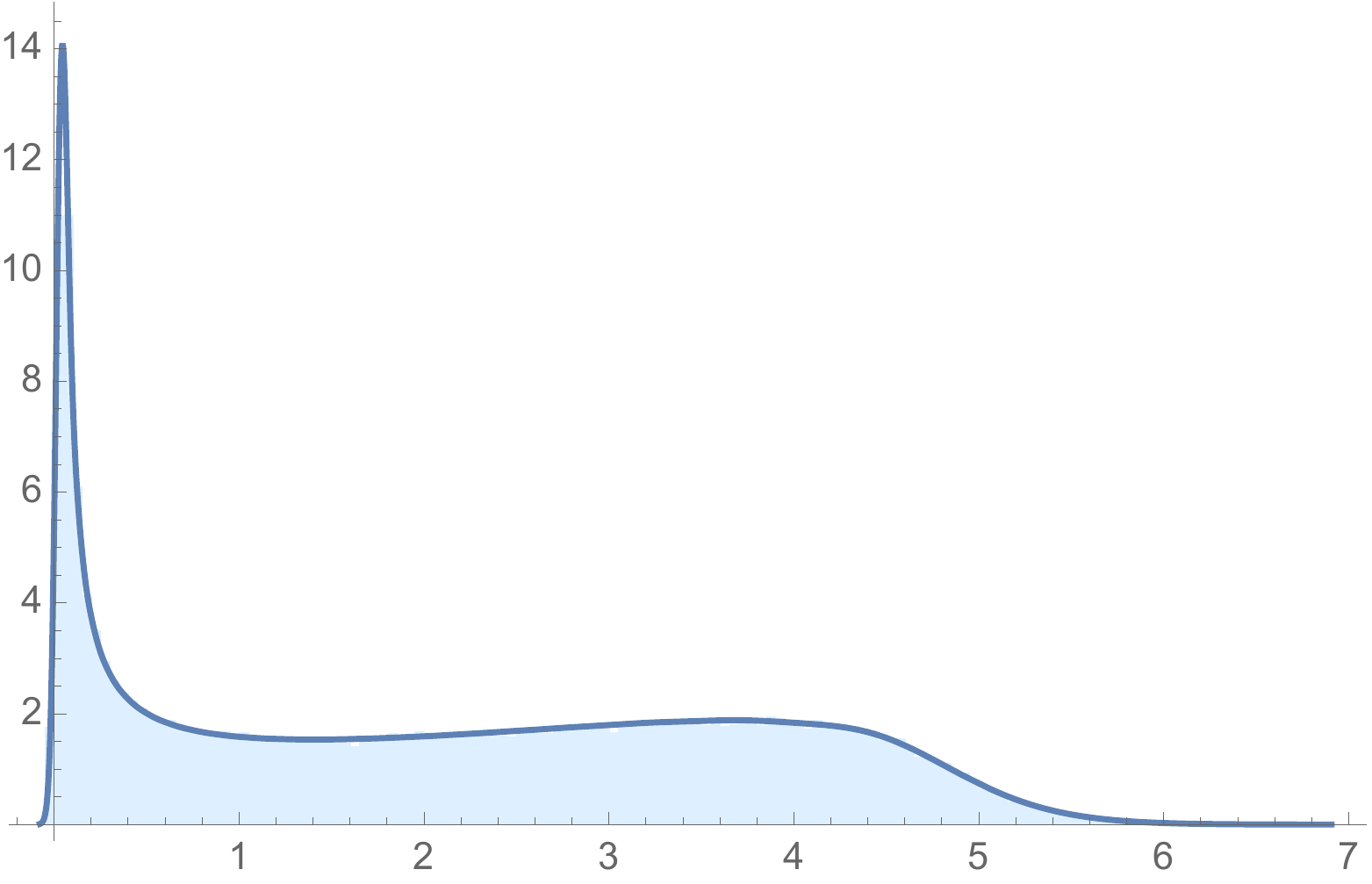}
\includegraphics[width=.30\textwidth]{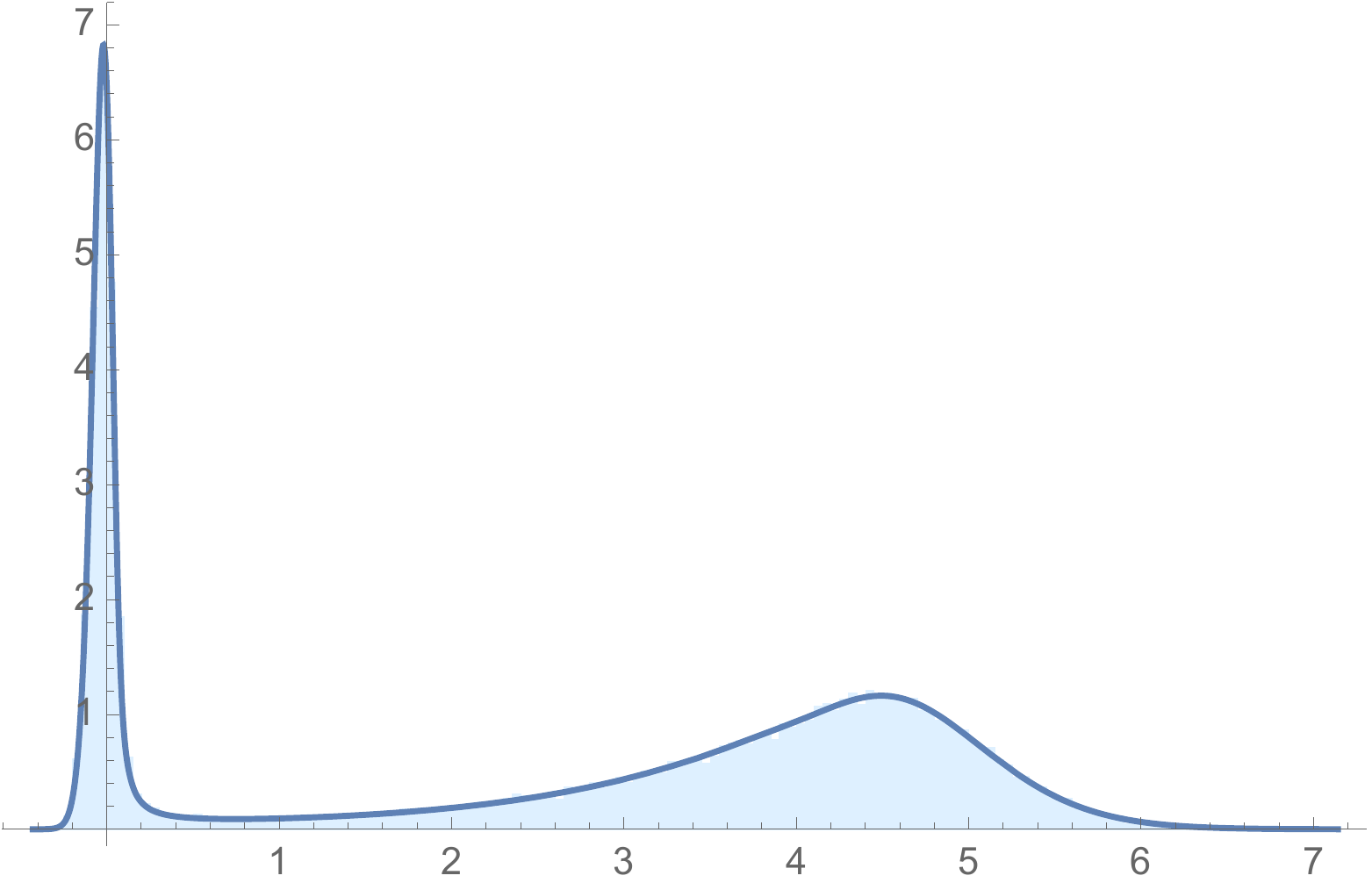}

\vspace{0.5cm}

\centerline{ Fig.2. Evolution of the density $p(\lambda)$ for $N=20\,,\,M=30$}
\centerline{ as the function of variance $\sigma^2=0.005; 0.25; 0.70$ }
\centerline{ {\small The blue histograms correspond to 10000 realizations.}}
\end{figure}

Our next goal is to investigate the limit $N$ $\&$ $M\, \rightarrow\infty$.

\subsection{ Asymptotic analysis.}

\subsubsection{\bf ``Bulk" Scaling Regime: extensive number of stationary points}

As $N$ $\&$ $M\, \rightarrow\infty$ in such a way that $1<\mu=M/N<\infty$ the number of stationary points in the loss function landscapes
shows three different regimes depending on the magnitude of the parameter $\delta=\frac{1}{2}\ln{(1+\sigma^2)}$.
The first regime is associated with the ``bulk scaling" corresponding to small enough
 $\delta\sim 1/N$ so that $\gamma = \frac{\delta N}{4}<\infty$. For such a regime one finds that the total number of solutions ${\cal N}$ is  {\it extensive}, namely:
 \be
 \lim_{N\rightarrow\infty} \frac{\mathbb{E}\{{\cal N} \}}{N} = \int\limits_{s_-}^{s_+} p_B(\lambda)\,d\lambda>0, \quad s_{\pm}=(\sqrt{\mu}\pm 1)^2,
 \ee
   where the density function  $p_B(\lambda)$ is expressed via  the  Marchenko-Pastur \cite{MP} limiting eigenvalue density $p_{MP}(\lambda)$ for the Wishart ensemble as
\be
  p_B(\lambda)= 2\, p_{MP}(\lambda)\,\exp{\left[-\frac{\gamma}{\lambda}(\lambda-s_-)(s_+-\lambda)\right]}, \quad
  p_{MP}(\lambda)=\frac{1}{2\pi \lambda} \sqrt{(\lambda-s_-)(s_+-\lambda)}
 \ee
 For $\gamma=0$ we obviously have $\mathbb{E}\{{\cal N} \}=2N$ whereas for $\gamma\gg 1$ we have asymptotically:
 \[ \lim_{N\rightarrow\infty} \frac{\mathbb{E}\{{\cal N} \}}{N} \Big|_{\gamma >>1} \approx \frac{1}{4\sqrt{\pi}} \frac{1}{\gamma^{3/2}}\ll 1
 \]
 Evaluating the above for $\gamma\sim N^{2/3}$ (i.e. $\delta\sim N^{-1/3}>>1/N$) indicates that the mean number of stationary points for such $\gamma$ becomes of order of unity as $N\gg 1$ defining a different scaling regime, cf. \cite{FLD2013}.

  \begin{figure}[h!]
\includegraphics[width=1.0\textwidth]{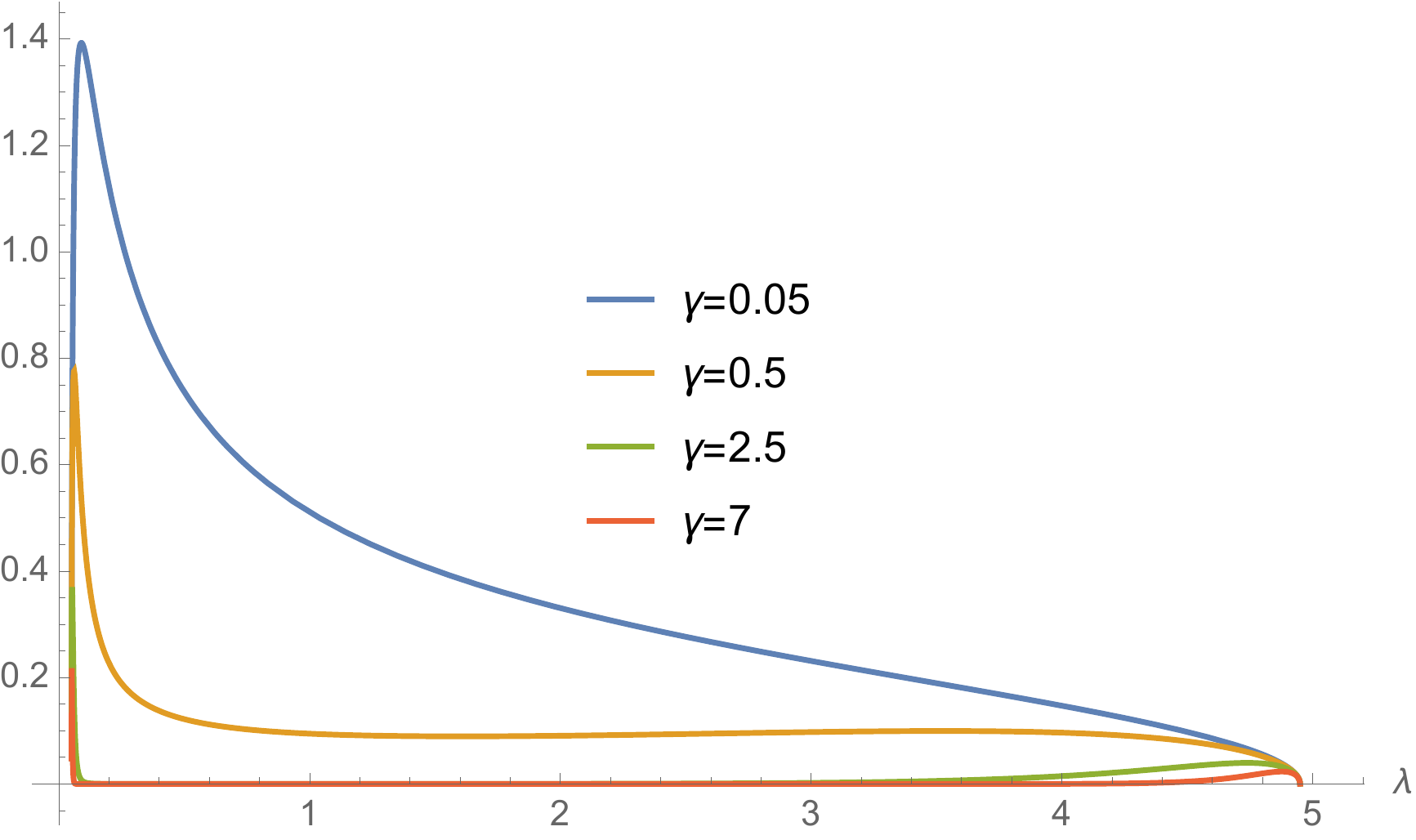}
\centerline{Fig. 3. Evolution of the density $p_B(\lambda)$ in the 'bulk scaling' regime.}
  \end{figure}

\subsubsection{\bf "Edge" Scaling Regime: finite number of stationary points}

The density of Lagrange multipliers for $\delta\sim N^{-1/3}$  is dominated by the vicinities of the spectral edges\\ \centerline{$|\lambda-s_\pm| \sim  N^{-2/3}\left(\frac{4s^2_{\pm}}{s_{+}-s_{-}}\right)^{1/3}\xi$}\\  where the Marchenko-Pastur law is no longer valid and has to be replaced by a more precise ``edge density" given by \cite{For12}
\be p_{MP}(\lambda)\longrightarrow \left(\frac{s_{+}-s_{-}}{4Ns^2_{\pm}}\right)^{1/3}\rho_{edge}(\xi),
\ee
with
\be
\rho_{edge}(\zeta)=\left[Ai'(\zeta)\right]^2-\zeta\left[Ai(\zeta)\right]^2+\frac{1}{2}Ai(\zeta)
 \left(1-\int_{\zeta}^{\infty}
 Ai(\eta)\,d\eta\right)
 \ee
where $Ai(\zeta)=\frac{1}{2\pi i}\int_{\Gamma}^{}e^{\frac{v^3}{3}-v\zeta}$ is the Airy function solving $Ai''(\zeta)-\zeta Ai(\zeta)=0$.

Introducing the scaling parameter $\omega=N^{1/3}\delta\left(\frac{s_{+}-s_{-}}{4}\right)$ one then finds the total number
of stationary points in this regime is finite as $N\to \infty$:
\be
\lim_{N\to \infty}\mathbb{E}\{{\cal N}\}=2\int_{-\infty}^{\infty} \left[\exp{\left(-\frac{\omega^3}{3s_{-}}+\frac{\omega\zeta}{s^{1/3}_{-}}\right)}+\exp{\left(-\frac{\omega^3}{3s_{+}}+
\frac{\omega\zeta}{s^{1/3}_{+}}\right)}\right] \rho_{edge}(\zeta)\,d\zeta
\ee
In particular, that number tends to just $\lim_{N\to \infty}\mathbb{E}\{{\cal N}\}=2$  as long as $\omega\to \infty$, indicating that for any
fixed and finite variance $0<\sigma^2<\infty$  only two stationary points typically exist:
one maximum and one minimum, cf. \cite{FLD2013}.

Comparison with results numerical simulations is shown in Fig.4.

\begin{figure}[h!]
\includegraphics[width=1.0\textwidth]{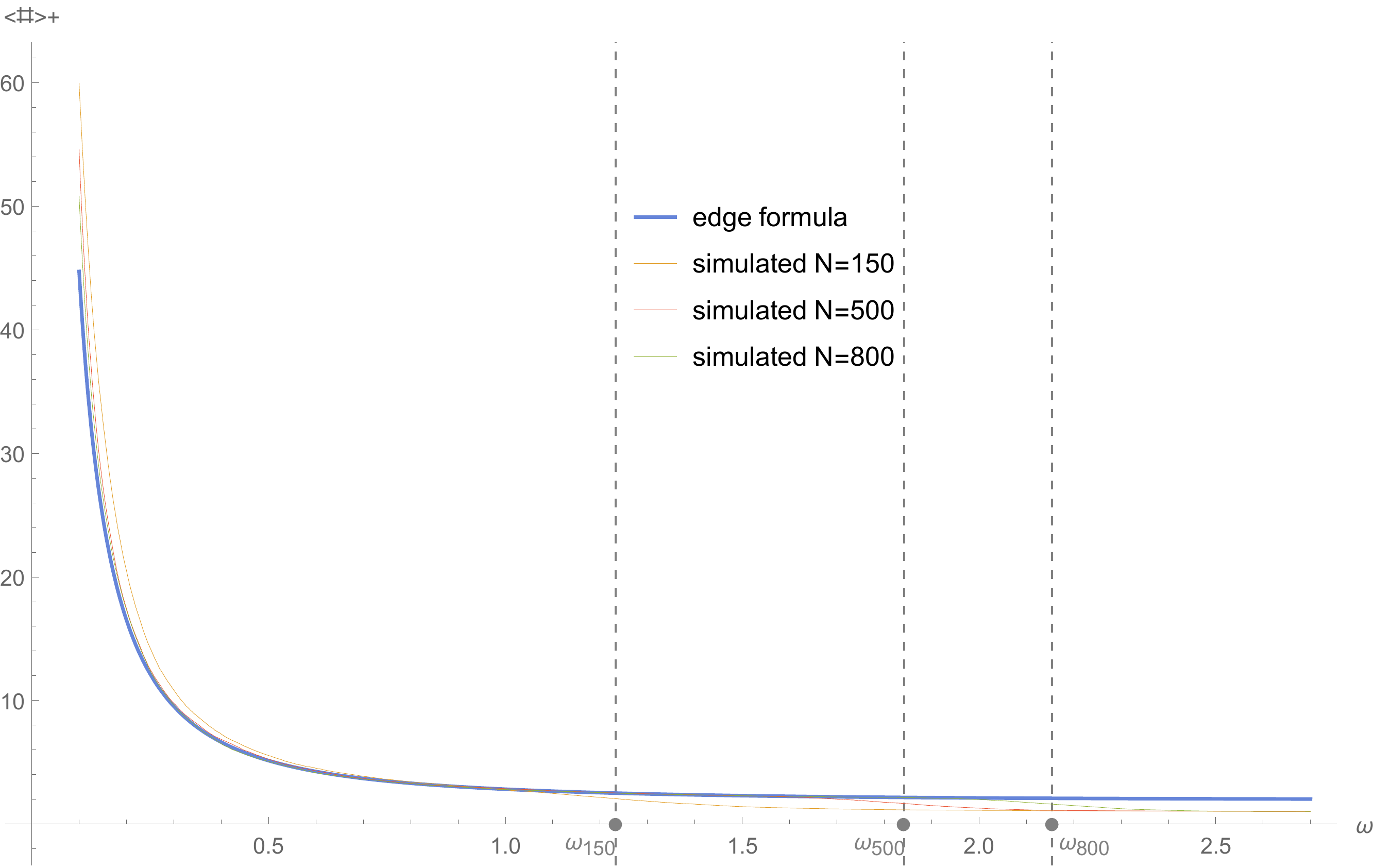}

\centerline{Fig. 4. Counting stationary points in the edge regime.}
\end{figure}

\subsection{\bf Large Deviations for the smallest Lagrange multiplier}

For large $ N\to \infty$, fixed $1<\mu=M/N<\infty$ and fixed finite $\sigma^2>0$ the probability density for the smallest Lagrange multiplier $\lambda_{min}$
has the {\it Large Deviation} form:
\be
p(\lambda_{min}<s_{-}) \sim e^{-\frac{N}{2}  \Phi(\lambda_{min})},
\quad \mathbf{ \Phi(\lambda)}= L_1(\lambda)+L_{2}(\lambda)+\frac{(\mu+1)}{2}\ln{(1+\sigma^2)},
\ee
where $s_{-}=(\sqrt{\mu}-1)^2$ is the 'Marchenko-Pastur' left edge and for $\kappa=\frac{(\mu-1)\sigma^2}{2\sqrt{1+\sigma^2}}$ we defined
 {\small \bea
 && L_1(\lambda)=(\mu-1)\left\{\frac{\sqrt{\lambda^2+\kappa^2}}{\kappa}-
 \ln{\left(\kappa+\sqrt{\lambda^2+\kappa^2}\right)}-
\lambda\frac{\sqrt{(\mu-1)^2+\kappa^2}}{(\mu-1)\kappa}\right\}\\
&& \mbox{and} \nonumber \\
&& L_2(\lambda)=-\sqrt{(\lambda-s_{-})(\lambda-s_{+})}
-2\ln{\frac{\left(\mu+1-\lambda+\sqrt{(\lambda-s_{-})(\lambda-s_{+})}\right)}{2\sqrt{\mu}}} \nonumber  \\
&&
+2(\mu-1)\ln{\frac{\left(\mu-1+\lambda+\sqrt{(\lambda-s_{-})(\lambda-s_{+})}\right)}{2\sqrt{\mu}}}
\eea
}

Comparison with the probability density of the smallest solution of Eq.(\ref{equlambda_new_xi_1}) found numerically is shown in Fig.5.

One then finds that $\Phi(\lambda)$ is minimized for
\be
\lambda=\lambda_*=(\sqrt{\mu}-\sqrt{1+\sigma^2})\left(\sqrt{\mu}-\frac{1}{\sqrt{1+\sigma^2}}\right)
\ee
providing the most probable/typical value of the smallest Lagrange multiplier. Substituting this value to
 Eq.(\ref{minLagr}) and then to Eq. (\ref{lossfun}) eventually gives the most probable value of the  minimal loss/error:
\be
\lim_{N\to \infty}\frac{{\cal E}_{min}}{N}=
\frac{1}{2}\left[\sqrt{\mu(1+\sigma^2)}-1\right]^2
\ee

\begin{figure}[h!]
\includegraphics[width=1.0\textwidth]{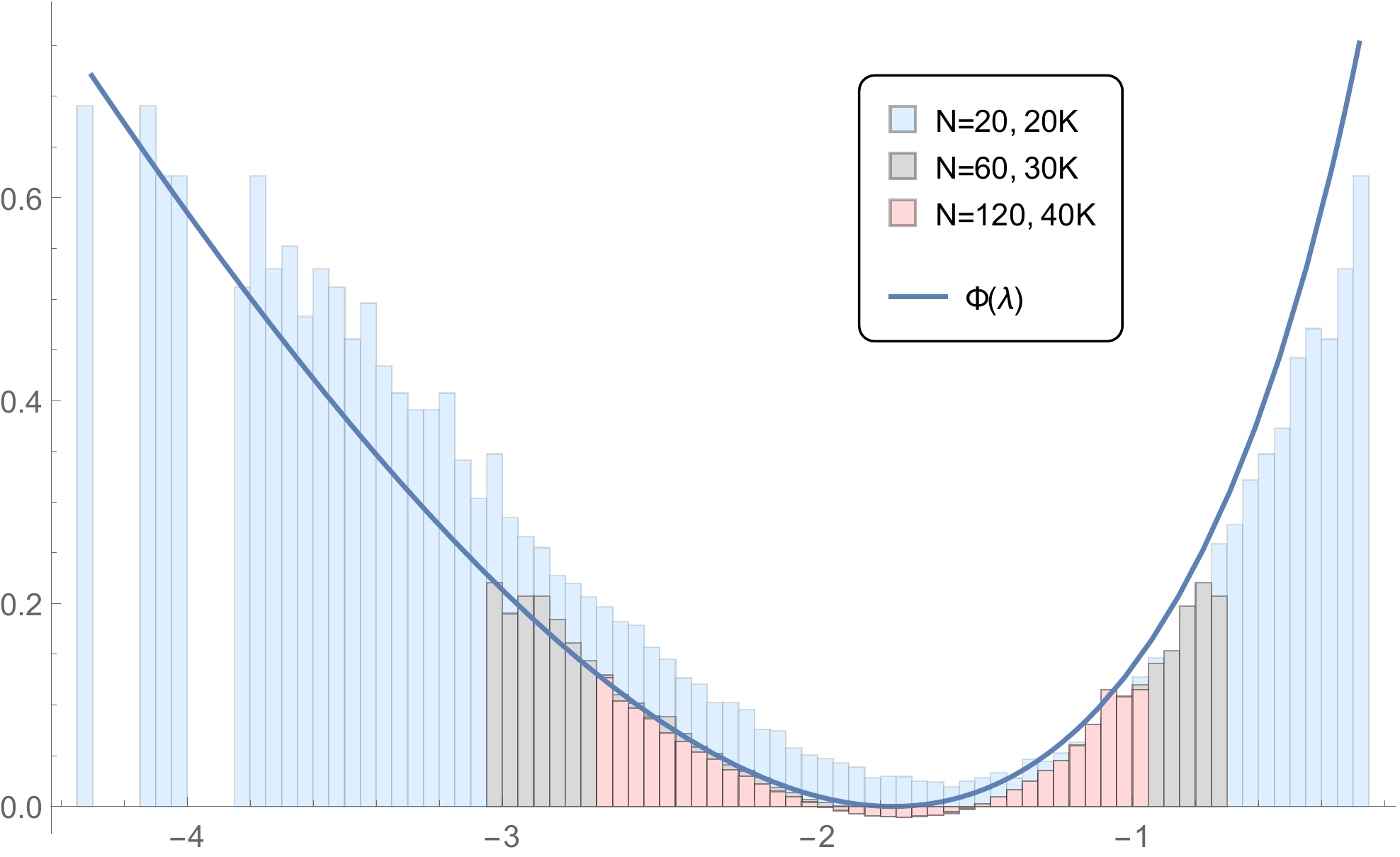}
\centerline{Fig.5. The large deviation function for the smallest Lagrange multiplier vs. simulations}
\centerline{ for different matrix sizes $N$ and different number of samples.}
\end{figure}

\subsection{Open questions}
 In conclusion, we counted the mean number of stationary points of the simplest ``least-square'' optimization problem on a sphere via  the  Lagrange multipliers in various scaling regimes, and  found the  typical minimal loss ${\cal E}_{min}$. The following questions remain open: (i) fluctuations of the counting function,
(ii) large/small deviations of the  minimal loss ${\cal E}_{min}$, (iii)  gradient search dynamics on the sphere
(iv) understanding the landscapes for  'least-square' optimization of more general type, e.g. involving nonlinearities etc,, cf. \cite{Fyo2019}. We hope to address some of these issues in future publications.

\end{document}